# Knot Logic – Logical Connection and Topological Connection
## by Louis H. Kauffman
## University of Illinois at Chicago

## I. Introduction

The purpose of this paper is to discuss how topology/geometry provides, in many instances, the connective tissue that enables logical comprehension. Probably the most well-known instance of this phenomenon is the Venn diagram approach to basic logic.

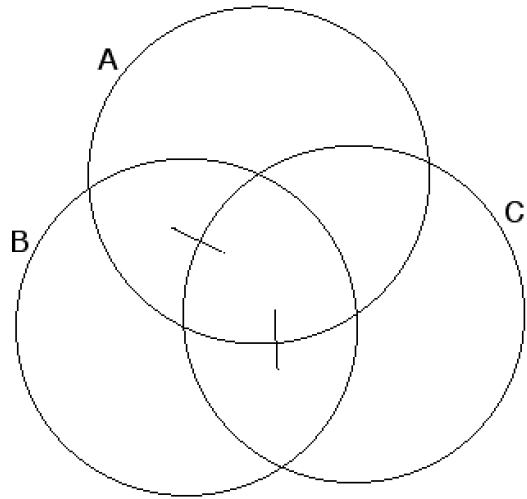

**Figure 1. Venn Diagram**

The diagram instantiates the eight possible states of affairs that can occur in a situation with three logical variables via the eight basic regions in the diagram. Statements about states of affairs regarding these variables can be indicated by marking the diagram. A collection of complex statements is recorded in the single diagram and the relationships among the statements can be read out from the marked diagram. The marked diagram circumvents step-by-step reasoning by its model of the very relationships about which one wishes to reason.

For example, in the diagram in Figure 1 we have markers in the regions where A intersects B and in the regions were B intersects C. Each marker spans two regions and is meant to indicate that *at least one of these two regions is not empty.* Thus the markers indicate the two statements

"Some A are B." and "Some B are C.".

It is then apparent at once that the statement "Some A are C." does *not* follow from these two statements since there is no necessity for the triple intersection of A,B and C to be occupied. In other words, one can easily read out a state of affairs that satisfies the first two statements, but not the third statement. All classical syllogisms and other matters of basic reasoning are resolved by the use of such diagrams.

We only indicate the use of Venn diagrams here. They are an exemplar of the sort of construction that interests us- *topological/geometric connection leading to the articulation of logical connection.*

The cartoon in Figure 2 illustrates another aspect of our theme. The little fellow (lets call him Parabel) in the cartoon is manipulating a unit interval via homotopy of smooth immersions, and he has discovered the "Whitney Trick" [1] that allows him to produce two opposite curls that move apart from one another and seem to have a life of their own. He has found topological entities (the two curls) that behave like +1 and -1, and he has shown us a topological version of the equation
$$0 = (-1) + (+1).$$
He has shown us that one can think of the production of +1 and -1 from 0 as a continuous operation that eventually gives birth to two independent discrete entities. Or we can think of the process (going backward) as a tale about how two independent discrete entities can interact continuously and cancel each other out. With these thoughts in mind, we are prepared to think further and wonder if the process, shown to us, might not be somehow analogous to the way a particle and an antiparticle interact physically to return to the vacuum!
In any case we see that there is a hint of an isomorphism of the integers with immersion homtopies of plane curves, and indeed this is the case. We are reminded by the topology, of the deep connections between the continuous and the discrete. And we see in this example how a continuous process can have discrete consequences via the discrimination of an observer. It is the topological observer who isolates the two curls and regards them as discrete entities separate from the continuous background.

Thus, when I speak of topological connection and logical connection, I do not just mean the diagrams or pictures, but also the way placing mathematical concerns in a topological context can show us larger connections that might otherwise have not been apparent .

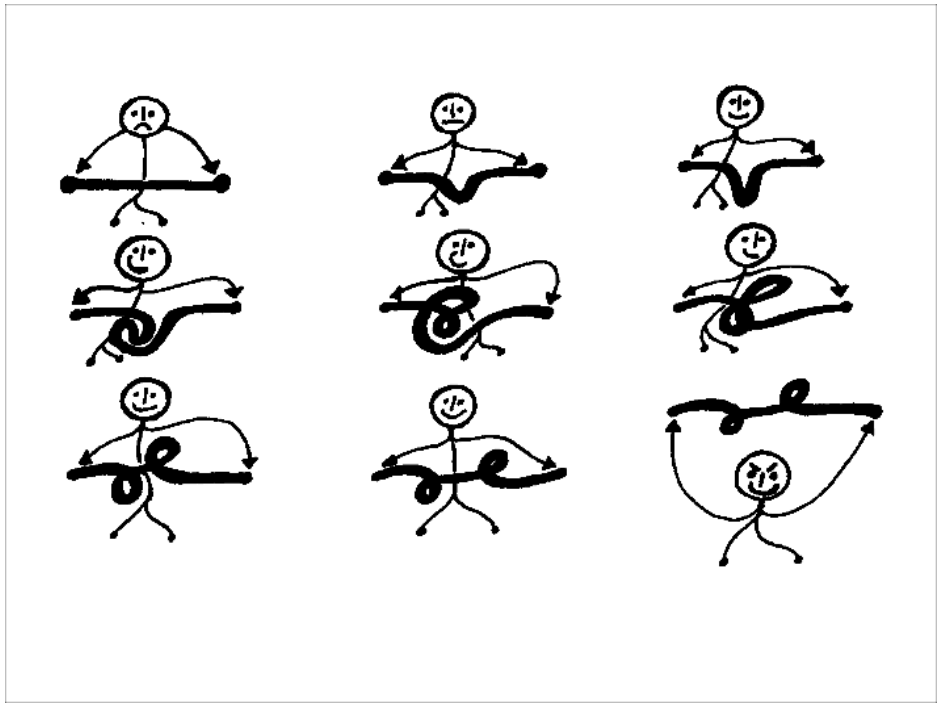

**Figure 2. Parabel making (-1) and (+1) from 0.**

**II. Knots and Knot Diagrams**

In Figure 3 we show a knot diagram. This should appear to the reader as a schematic sketch of a knotted closed curve. But it can also be seen more formally as depicting a plane graph with order 4 nodes that have been each given extra structure in the form of a "crossing", where two out of the four local edges at the node are connected and the other two are disconnected by a broken line segment. The two edges that are connected appear at two places of the same parity in the planar cyclic order of the node. Again, pictorially, the crossing has the appearance of one part of the knot crossing over another part. The diagram has a multiple use. It is a sketch of the knot, and can be used to weave an actual knot from rope corresponding to this sketch. It is a depiction of a formal graphical structure that can be taken to represent the knot mathematically. These two fundamental aspects of the diagram are equally important. They serve to connect the abstract structure of the graphical model with the corresponding abstract structure of the knot defined to be an embedded curve in three dimensional space. They also serve to connect the abstract structure with a recipe for making a knot model from rope, or for recording a knot originally made in rope into the mathematics. The graphical structure can be translated to provide computer programs with the information to compute invariants of the knot.

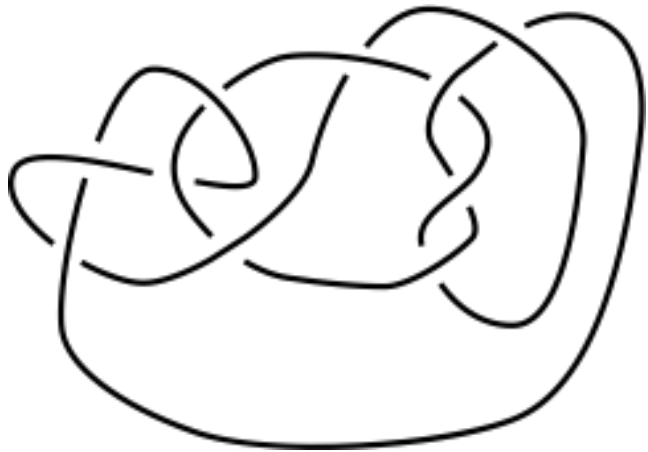

**Figure 3 - A knot diagram**

In Figure 4 we illustrate the *Reidemeister Moves*. These are moves that can be performed locally on knot or link diagrams that preserve the topological type of the corresponding curve in three dimensional space. From the point of view of the graphs, move I removes or adds a 1-sided region, move II removes or adds a 2-sided region and move III changes the relationships of the diagram by replacing a 3-sided region by another in the course of rearranging three crossings. The moves are performed with the crossing relations as shown (plus some versions of the same type that we have not illustrated here). In the 1920's Kurt Reidemeister [2] proved that two knots or links are topologically equivalent in three dimensional space if and only if their diagrams can be transformed to each other by a sequence of the Reidemeister moves. These moves complete the connection between the abstract graph model for knot theory and the topological concept of knot equivalence by isotopy in three dimensional space.

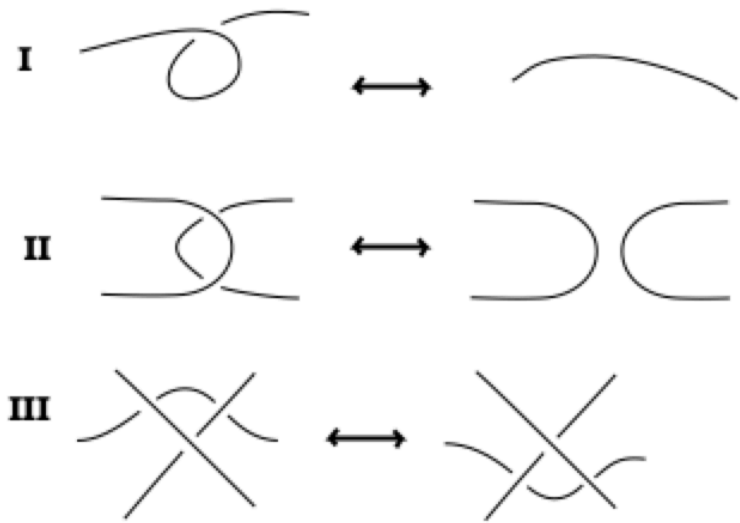

**Figure 4 – The Reidemeister Moves for Knot and Link Diagrams**

As we will see later, the value of the Reidemeister moves is theoretical in that their effects can be checked easily due to their simplicity. It is for that reason that they have been effective in providing a connection between the topology of knots and links and other subjects such as algebra, group theory, combinatorics and statistical physics. In an actual practical problem about undoing a knot, it is helps to use larger scale moves. We will illustrate such moves momentarily. First consider the problem indicated in Figure 5. Is the "knot" shown in Figure 5 actually knotted?
In Figure 6 we show a diagram that has the same weaving as the three dimensional image in Figure 5. The unknotting is shown in Figure 7 with the help of knot diagrams and a large-scale move-type that I call a "swing move".   The first frame in Figure 7 illustrates the same diagram as Figure 6, with two points marked by dark dots.   An arc between the two dots goes entirely under the rest of the diagram, and it can be "swung" underneath the diagram to the position in the second frame in Figure 7 where it only undercrosses with a single crossing. *The swing move itself can be accomplished by a series of Reidemeister moves. We leave this factorization as an exercise for the reader.* In the third frame we indicate another arc that goes over three times, and can be swung to eliminate the three crossings.  After that, two type II moves reduce the diagram to the next-to-last frame and this is transformed to the unknot by a move of type I (and some smoothing out of the resulting diagram). Thus the knot in Figure 5 is seen to be unknotted.

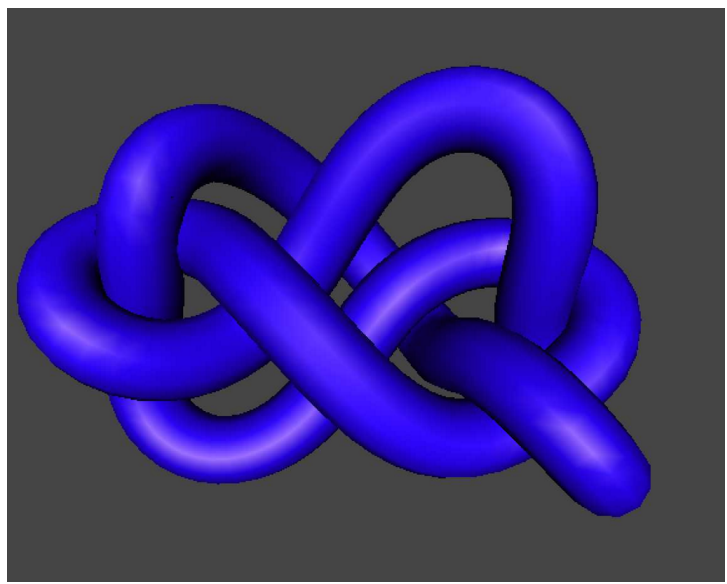

**Figure 5 – Is this knot knotted?**

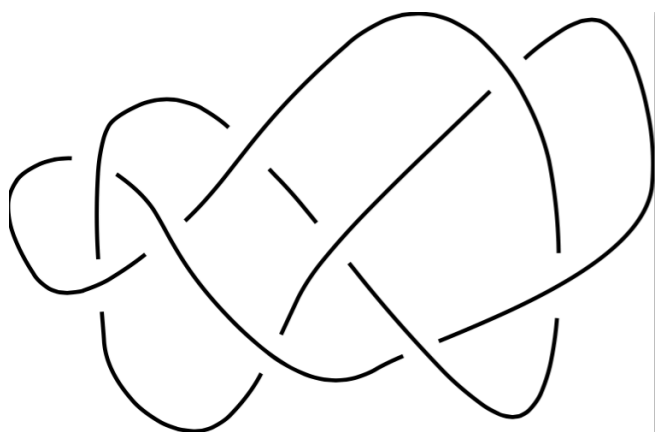

**Figure 6 – A knot diagram for the knot in Figure 5.**

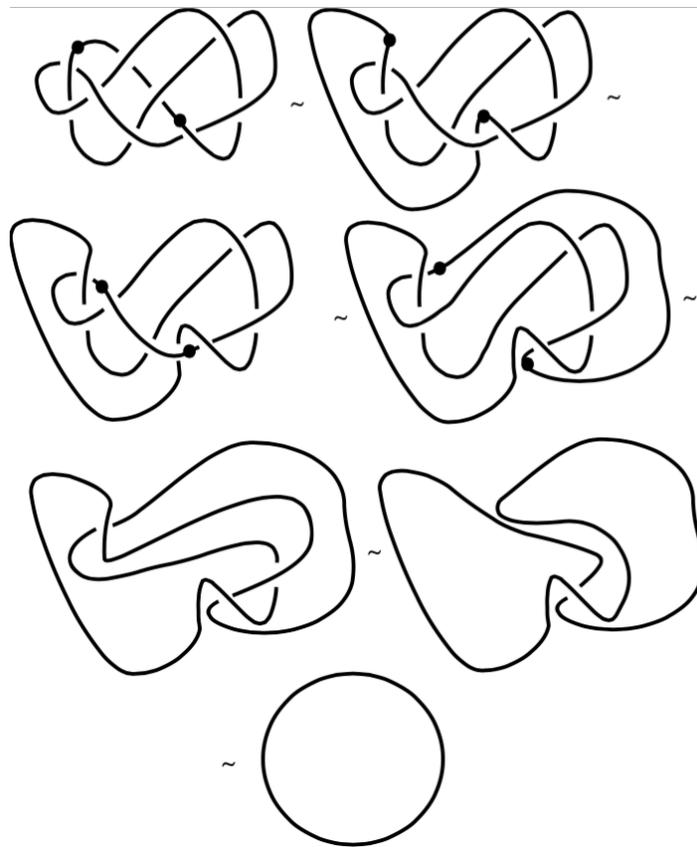

**Figure 7 – The unknotting of the diagram in Figure 6.**

We see from this example, that the (graphical) knot diagram, coupled with topological moves on diagrams, provides a language in which a person can record and indeed discover significant relationships about knots. The diagram is an intermediary between an abstract mathematical model of knotting and the actuality of knotted objects in three dimensional space.

In Figure 8 we show an electron micrograph of a closed circular DNA molecule that has been coated with protein [3] in such a way that the knotting of the DNA is apparent from the micrograph. In fact, the electron micrograph has been made in such a way that a biologist can read it as a knot diagram. The concept and application of knot diagrams was essential for this scientific application. Electron micrographs are two dimensional renderings. It was necessary that topological information in the DNA knot could be determined from such a projection . The knot diagram becomes a connection between the invisible world of the DNA molecules and the structural topological world of the mathematics.

Now examine Figure 9. In this figure we begin with an idealized bit of closed circular DNA, not knotted, but harboring a three-fold twist. Just after the arrow, the DNA is bent over so that two of its arcs are in proximity, and its knot type is indicated as a closed circle to the right. At the second arrow a recombination event occurs. The DNA is broken and re-spliced, forming a crossing of right-handed type. The result of the recombination is a simple link of two components with linking number equal to 1. A second recombination occurs and a figure-eight knot appears. Then comes a so-called Whitehead link and then a more complex knot. We see in this diagrammatic experiment the results of successive recombination, where the pattern of recombination (insertion of one right-handed twist) is always the same. The diagrammatic experiment is linked (no pun intended) with real experiments and appropriate electron micrographs [3] to show that a real processive recombination actually does produce these knots and links. The result is that the hypothesis of the form of the recombination is confirmed by this diagrammatic experiment coupled with the surrounding knot theory and the surrounding molecular biology. It is the language of the diagrams that provides the crucial connection between the biology and the mathematics.

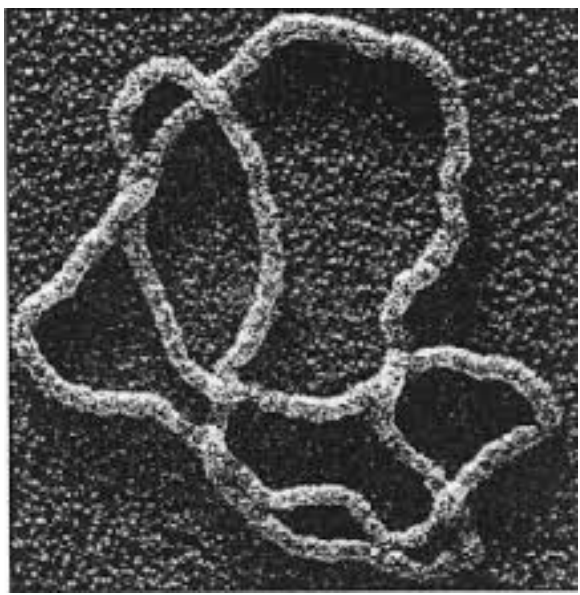

**Figure 8 – An electron micrograph of knotted DNA.**

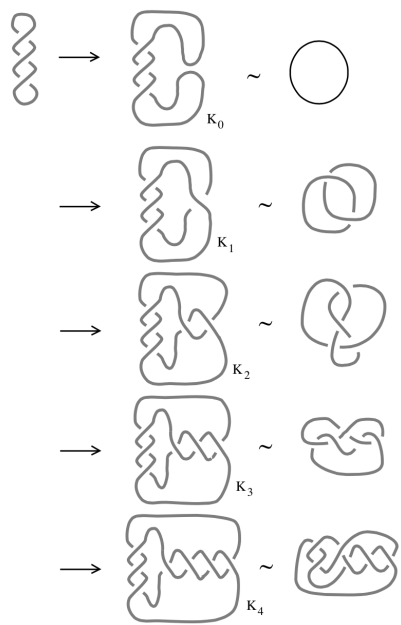

**Figure 9 – Diagrammatic Experiment in Processive DNA Recombination**

### III. The Dirac String Trick

The Dirac string trick is a property of braiding that can be illustrated with a twisted belt. The application is to the quantum mechanical properties of an electron (or other Fermi particle). The trick with the belt is illustrated in Figure 10.

## The Dirac String Trick

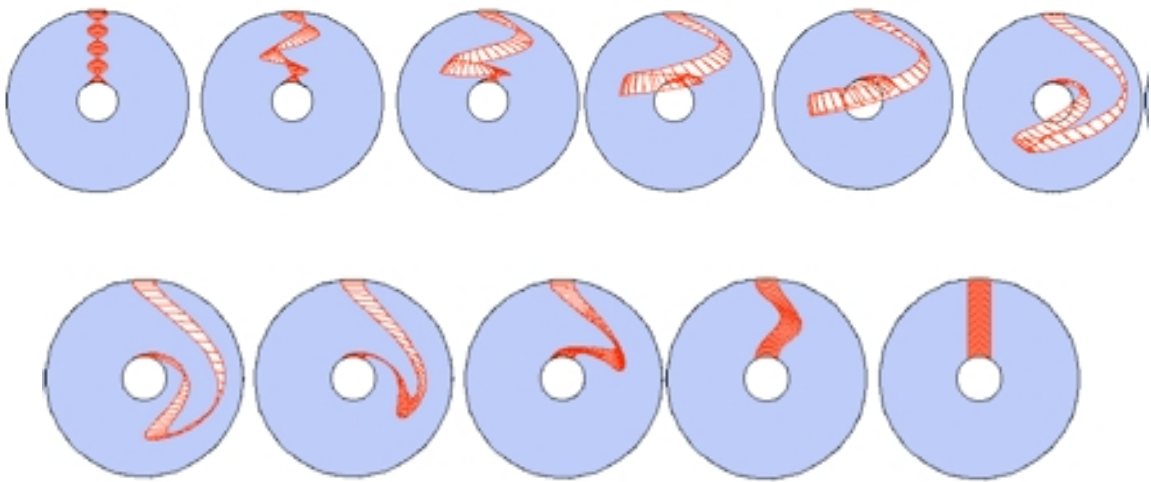

Figure 10 – The Dirac String Trick

In this figure we see a belt attached between two concentric spheres. Throughout the indicated motions, the belt is anchored to both the inner and outer spheres, but it is free to move elastically around the inner sphere. By sweeping it around the inner sphere, the belt returns to its original position with all of the twist removed.

This topological phenomenon with the belt can be demonstrated quite easily with a real belt. From a mathematical point of view it can be seen as a visualization of the fact that the fundamental group (in the sense of algebraic topology) of SO(3), the group of axial rotations of three dimensional space, is isomorphic to Z/2Z, the integers modulo two. This, in turn comes from the fact that SO(3) is topologically equivalent to the three dimensional projective space. To see this fact, view Figure 11.

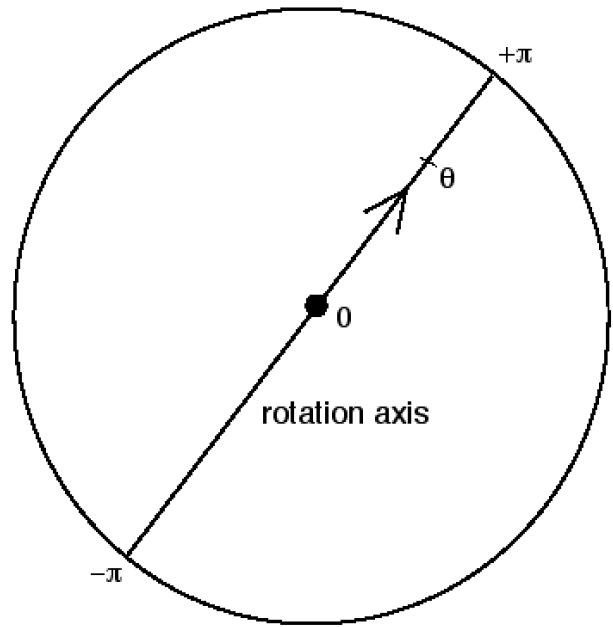

**Figure 11 – SO(3) is topologically the three dimensional projective space.**

In Figure 11 the angle of rotation is marked on the rotation axis, from -π to +π on a ball of radius π. Thus SO(3) is topologically a ball with antipodal points on its boundary identified with one another. The identification space is homeomorphic to the three dimensional sphere with antipodal identifications. This is the three dimensional projective space.

Note that each diameter of the ball corresponds to a loop of rotations in the rotation group by using the angle parametrization of the points on that diameter. An individual configuration of the belt can be obtained by starting with a straight vertical belt between the two spheres and taking its image under a loop of rotations that starts at the identity rotation and ends at the identity rotation. In this way the first image in Figure 10 is the result of the square of the generator of the fundamental group of SO(3), and the other images represent a null-homotopy of this element of the fundamental group. For more information about this aspect of the belt trick and the related physics, the reader should see [4].

The physicist Paul Dirac liked to use the belt trick to illustrate a strange phenomenon about electrons. If you circulate by $2\pi$ around an electron, the quantum state of the electron changes by a phase factor of (-1). In other words, from the point of view of an electron, a $2\pi$ rotation is different from no rotation at all. But a $4\pi$ rotation comes back to the identity! The reasons for this strange behavior have to do with the fact that quantum processes are represented by unitary transformations and the corresponding unitary group for the electron is SU(2) whose topology is the three-dimensional sphere. SU(2) is the double cover of the rotation group SO(3). The result is that the topological and group theoretical mathematics of the belt trick is directly relevant to the quantum states of an electron. It is as though the electron were connected by invisible strings to its surroundings in such a way that these strings unwind themselves after every $4\pi$ rotation.

The Dirac string trick is a deep example of the way topology can participate in both topological and logical connection of apparently disparate subjects. There is much more to say about this subject. See [4]. A few more words are appropriate here. The group SU(2) is isomorphic with the group of unit quaternions, of the form a + bi + cj + dk where {1,i,j,k} are Hamilton's quaternion generators satisfying ii = jj = kk = ijk -1. See Figure 12.

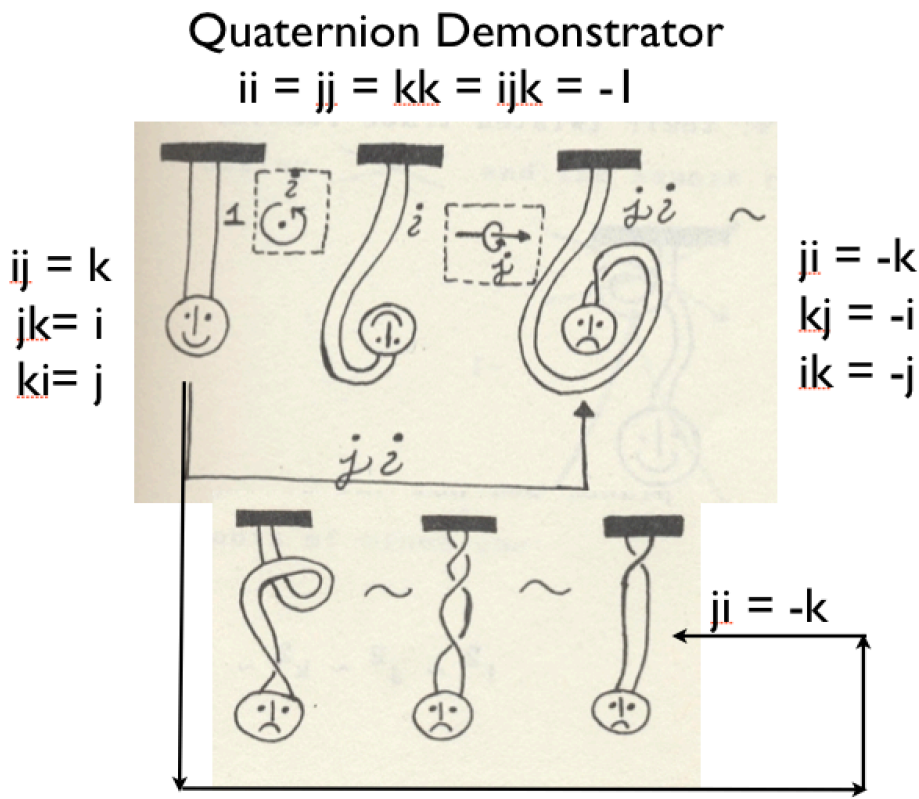

**Figure 12 – The Quaternion Demonstrator**

In Figure 12 we illustrate how the belt trick is related to the quaternions.
A twist of π on the belt can be done around any spatial axis. Let i,,j and k denote such twists of π around three perpendicular axes. Now the belt is attached to a wall at one end and a card (with cartoon faces on it in the figure) at the other end. Rotate the card around a given axis. The π rotation is always a topological square root of minus one, since four such rotations return the belt (with a little topological help) to its original position. Then one finds that the composition of these rotations, combined with the state of the belt, give rise exactly to the quaternion relations! You can do it with a belt. You can do it with your hand and arm. You have the quaternions in the palm of your hand.

**IV. Knot Diagrams and Electricity**

In this section we describe a surprising connection between knot diagrams and electrical networks. View Figure 13 and you will see a translation of the Reidemeister moves to moves on a graph that is associated with the knot or link diagram.

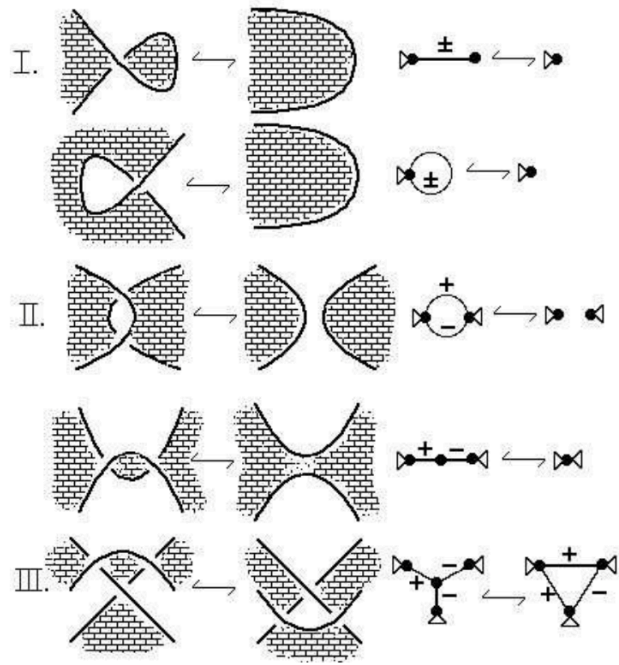

**Figure 13 – Graphical Translation of the Reidemeister Moves**

The associated graph is obtained by shading the knot diagram in two colors so that adjacent regions are colored differently. Each region of the checkerboard coloring becomes a node of the associated graph. Two nodes are connected when there is a crossing shared by the two regions to which they correspond. An edge of the graph is labeled with a +1 or a -1 according to the convention implicitly illustrated in Figure 13. The Reidemeister moves then translate into the addition or removal of a pendant edge or loop (for move I), series and parallel transformations (for move II) and a star-triangle relation (for move III). With the patterns of signs as shown in Figure 13, *these are moves that preserve the conductance of the associated graph viewed as an electrical network.* See Figure 14 for a summary of the more general rules for handling equivalent conductance in electrical networks.

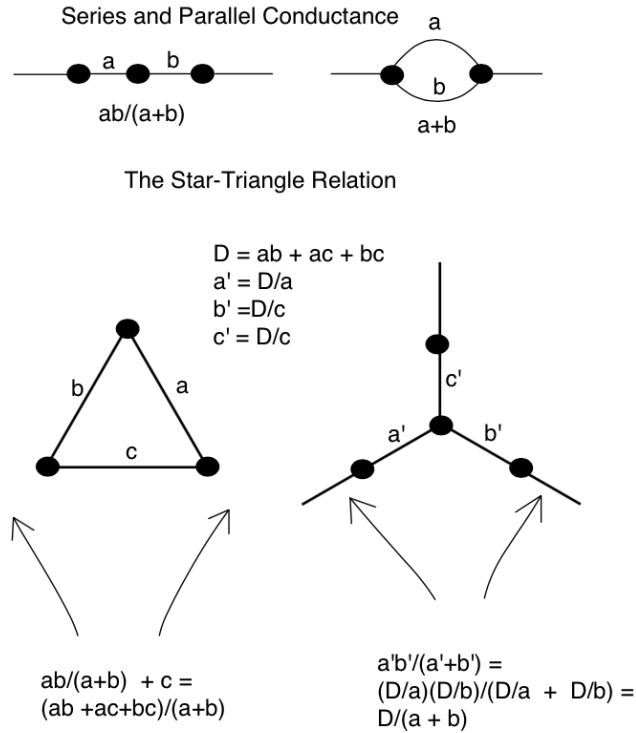

**Figure 14 – Series, Parallel and Star-Triangle Moves for Electrical Networks**

For example, consider the Borommean rings shown in Figure 15. These rings are linked as a triple, but any two of them (in the absence of the third) are unlinked. We have drawn the corresponding electrical graph, and see at once that the conductance will be non-zero for any two nodes as terminals, since all the edges have conductance +1. This means that we have proved that the Borommean rings cannot be unlinked by any moves that do not pass the diagram across the chosen terminals.

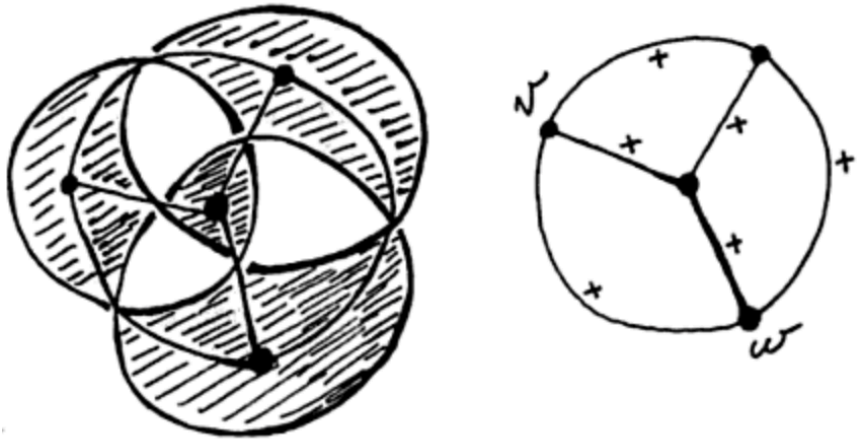

**Figure 15 – Conductance Network for the Borommean Rings.**

This remarkable correspondence of knot diagrams and electrical networks can be used to obtain topological information about knots and it can be used to analyze networks. This correspondence, discovered by L. Kauffman and J. Goldman [5], has not yet been fully understood. It is a topological connection not all of whose connections are yet apparent.

**V. Knot Logic**
In this section, we discuss very general themes related to knotting and linking. In the next section we give a formal model for a non-standard set theory, using knot diagrams, that allows self-membership and mutual membership. We begin with Figure 16, illustrating the concept of linking as a mutual relationship among its components.

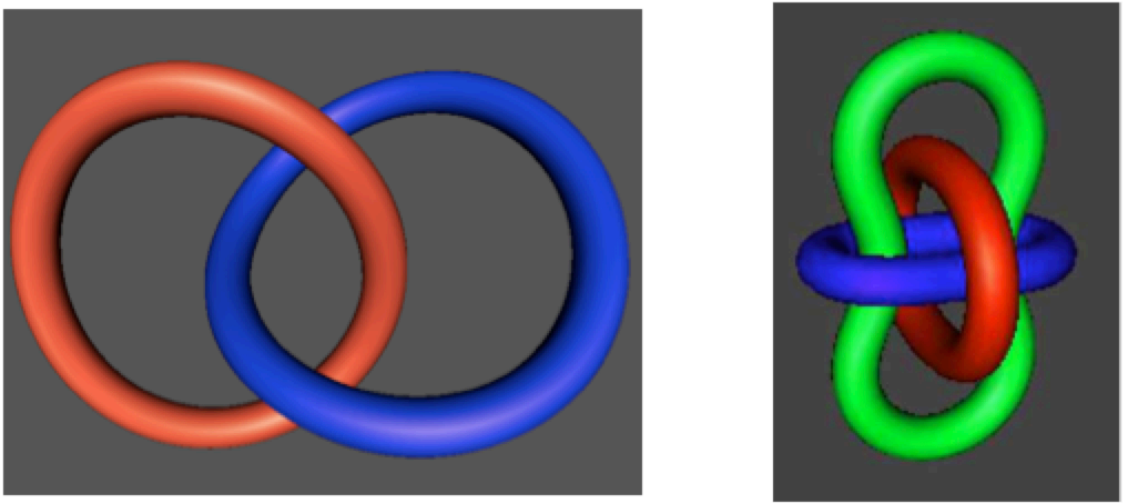

Figure 16 – Linking as Mutuality

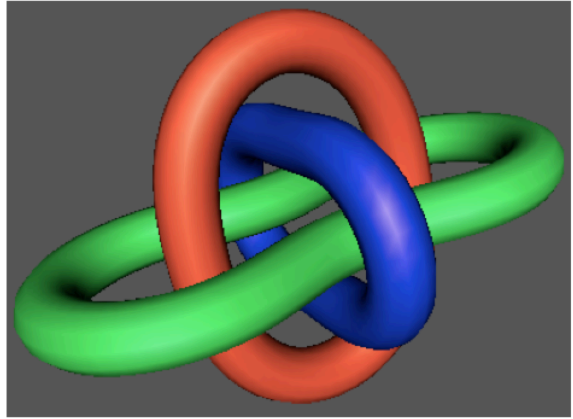

Green surrounds Red.
Red surrounds Blue.
Blue surrounds Green.

**Figure 17 – Borommean Rings and The Surround Relation**

In Figure 17 we show a three dimensional image of the Borommean Rings in which it is apparent that any two rings are unlinked. The rings are colored Red, Blue and Green. One can see that Red surrounds Blue, Blue surrounds Green and Green surrounds Red! This cyclicity in the surround relation is certainly related to the fact that the rings themselves are linked. But at this level, we are noticing facts about this particular view of the Borommean rings. It is not so easy to make a version of the knot theory where this notion of "surrounding' can be given topological meaning.

Figure 18 carries an idea that was originally made popular by Buckminster Fuller [6] that he called *pattern integrity*. One says that a knot has pattern integrity to mean that *as a knot* it is independent of the structure or composition of the substrate that carries it.

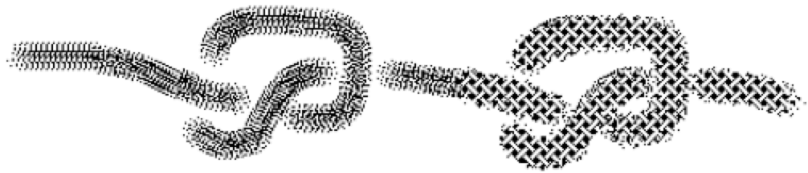

**Figure 18 – Pattern Integrity – Independence of Substrate**

This idea is an important one for thinking about systems that are supported by a substrate, such as biological systems. Such systems have a pattern integrity of their own that transcends the components and parts that are continually changing in the course of time. It is the pattern that persists and it is the tendency to preserve certain patterns that gives the organism its integrity. We realize on looking at this concept from architecture, biology, systems theory and cybernetics that it is at the core of our topological considerations. We are not only concerned with properties of the knot that are independent of the substrate, we are also concerned with independence of the shape of the substrate up to certain transformations (isotopy or the Reidemeister moves) and so we are engaged in a theory of forms, of pattern integrity that can apply to systems and to biology exactly through the connections implicit in topology and geometry.

## Self-Mutuality and Fundamental Triplicity

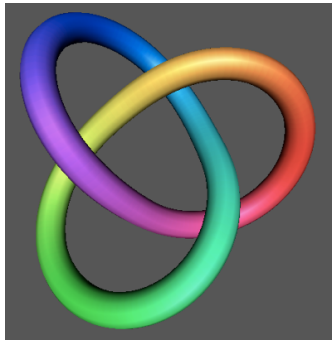

Trefoil as stable self-mutuality
in three loops about itself.

Figure 19 – Trefoil as Stable Self-Mutuality – Interlocking Its Own Loops

View Figure 19. When you weave a trefoil knot, you make a loop and then you pass the arc of the knot through that loop. The arc passed through stabilizes the original loop and makes the knot stable except for the possibility of sliding off the ends of the rope. Then with the self-reference of splicing the rope into itself to produce a closed knotted loop, you have the knot fully independent of your constructing hands, free to remain knotted in three space. The trefoil knot is a stable self-mutuality in three loops about itself.

## VI. Knot Set Theory

In this section we create a non-standard set theory using the formalism of knot and link diagrams. We begin by viewing crossing as a membership relation. After all, if b crosses over a then this is not symmetric, and membership is a non-symmetric relation. See Figure 20 for the basic set-up. We label the arcs of the diagram as shown and say that an arc *a belongs to* an arc b if a under-crosses b.

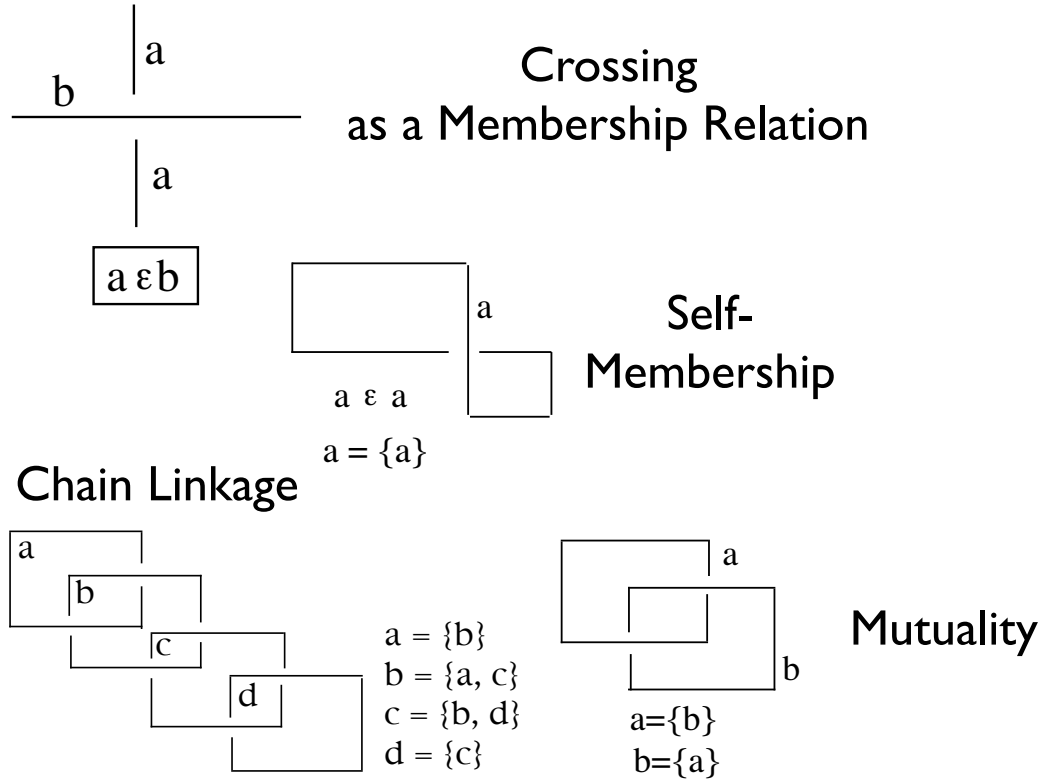

**Figure 20 – Crossing as a Membership Relation**

Figure 20 shows how non-standard set theoretic relations are quite natural from this point of view. In the figure you will see a diagram for self-membership a = {a} as a curl, and mutuality a = {b}, b = {a} as linking. A chain linkage has the form a = {b}, b = {a,c}, c = {b,d}, d = {c}.

It is clear that we can use the knot and link diagrams to illustrate these iconoclastic notions. But what about the topology of the knots and links? View Figure 21. We see that the second Reidemeister move demands that we cancel repeated elements of a knot-set in pairs. Thus the set {a,a,a} will be equal to {a} as is usual in multi-set to set conversion, but {a,a} will be equal to the empty set { }. We will therefore take our knots sets under this Fermionic reduction rule. Their diagrams will then be invariant under the second Reidemeister move. There is no problem with the third Reidemeister move. The membership relations remain the same.

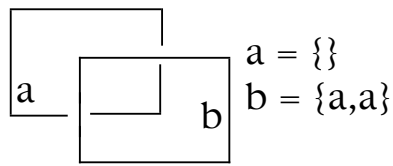

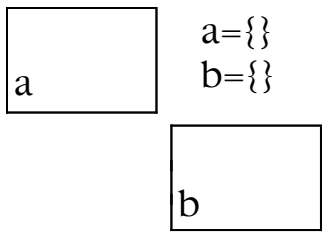

**Figure 21 – The "Fermionic Nature" of Knot-Sets**

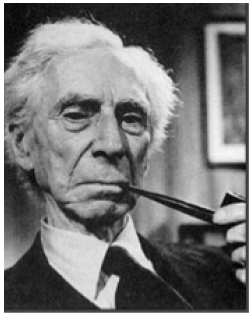
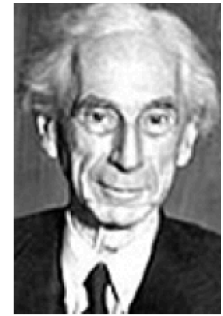
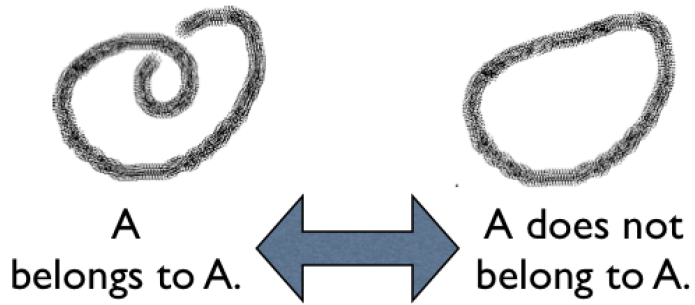

**Figure 22 – The Russell Paradox and the First Reidemeister Move**

Finally, in Figure 22, we illustrate the situation with the first Reidemeister move. If we allow it, then it says that a set that is its own member can retain this membership or have it removed. The two situations are equivalent. It would seem that the Russell paradox vanishes, since then all sets are or are not members of themselves as you wish! A closer analysis of the paradox reveals that we shall still have to consider the paradoxical nature of the set of all sets. This analysis will be published elsewhere. The reader who examines our figure closely will see that we have written the equation Rx = ~xx. This is shorthand for "x is a member of R if and only if x is not a member of x." R is the Russell set.

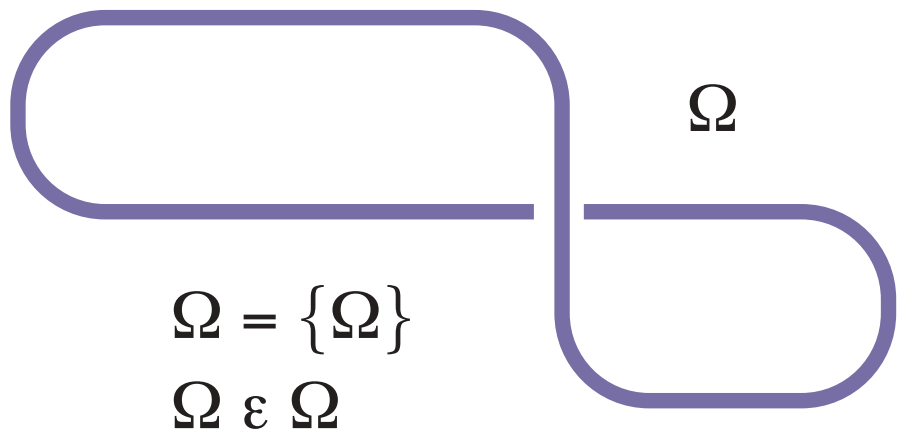

**Figure 23 – The Set whose only member is itself.**

Another option for knot-sets is to allow only the second and third Reidemeister moves. In knot theory the equivalence relation generated by moves II and III is called regular isotopy. We shall call knot-sets under these two moves *framed knot-sets*. Figure 23 illustrates the basic framed knot-set that is a member of itself and has no other members. The figure shows how we can think of
$$\{\Omega\}=\Omega$$
without infinite regress. One can imagine an observer walking around the diagram and changing role from member to container to member again and again throughout the circular walk.

A further comment:. Knot sets can be shown as tangle diagrams and if there is no weaving, then the result is a classical set. For example in Figure 24 we have illustrated the recursive construction of the ordinals due to John von Neumann.

$$0 = \{\ \}$$
$$1= \{0\}$$
$$2=\{0,1\}$$

3={0,1,2}

Architecture of Counting

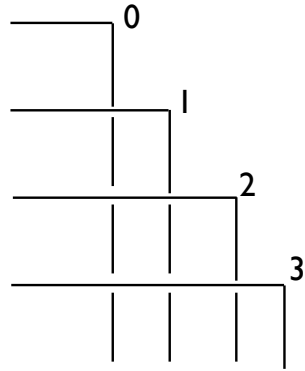

**Figure 24 – Ordinals constructed recursively.**

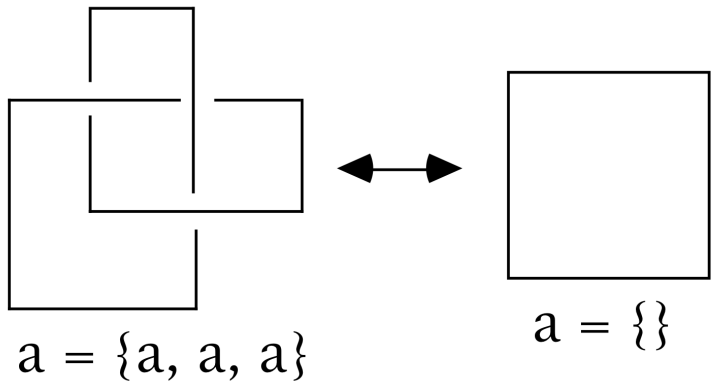

**Figure 25 – The Trefoil Knot-Set is Unknotted.**

Returning now to knots, we realize that topological knot-set theory
is in fact only applicable to links. See Figure 25. Due to the Fermionic nature of the
theory, once we allow all three Reidemeister moves all knots will disappear.
Even if we only use moves II and III all knots will still disappear, leaving only some
curly patterns of self-membership.

Thus in order to see knots we must use a stronger form of logic.
And it is natural from this point of view to wonder if in our logical explorations we
have missed some form of logical thinking that is more closely related to topology.
That is the subject of the next section.

## VII. Discrimination Logic, Quandles and the Trefoil Knot

Imagine a very elementary universe in which there are two distinct entities A and B that can interact with themselves and with each other. Suppose that A interacts with itself to produce only itself. We write AA = A. And suppose that B interacts with herself to produce B. We write BB = B. Now if A and B were to interact and produce Nothing we would have a very elementary universe indeed with A and B antiparticles. Assuming that A and B existed, it would not be long before there was only Nothing! So let us suppose that A and B can interact to form a new entity C.
AB = BA = C. Now it is possible to have many new entities if, for example, A were to interact with C to produce a new entity D and so on. We are familiar with the generation of infinitely many new sets from nothing but the empty set (and the act of collection). But here we choose simplicity. Let A interact with C to produce B and B interact with C to produce A. Then we have
$$AB = BA = C,$$
$$AC = CA = B,$$
$$BC = CB = A.$$
There are just three elements in this universe. Each pair of elements maintains the distinction that is the pair by interacting to produce a third entity. Conceptually, this 3-Universe is as simple as it can be and still not collapse to Nothing. See Figure 26.

We can certainly agree that this *3-color algebra* is fundamental, primordial, at the very beginning of our mathematical considerations. And yet notice this
$$A(BC) = A(A) = AA = A,$$
$$(AB)C = C(C) = CC = C.$$
Thus A(BC) is not equal to (AB)C and the 3-color algebra is not associative.

A Logic of Elementary Discrimination.

Entites A and B are distinct.
They interact to produce a new entity C.
A is indistinct from itself and interacts
with itself to produce itself.
AA = A, BB= B, CC=C.
The system of interactions closes.
BA = AB = C
CA = AC=B
BC = BC = A

**Figure 26- A Logic of Elementary Discrimination**

Remarkably, this 3-color algebra is directly related to knots and to the Reidemeister moves. View Figures 27 and 28. In Figure 27 will illustrate how to locally color the edges of a knot or link diagram using the algebra. An undercrossing line labeled X meets an overcrossing line labeled Y and continues on the other side of the overcrossing line labeled with the product XY.

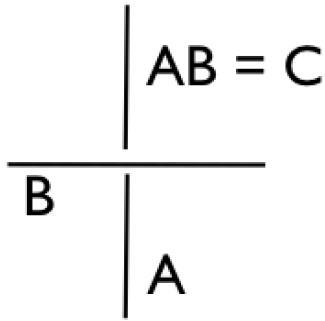

Now the arcs in the diagram are colored by A,B,C.

The algebra of interactions of A,B,C is not associative! It should remind you a bit of the vector cross product.

$$(AA)B = AB = C$$
$$A(AB) = AC = B$$

Figure 27 – Three Color Logic on the Knot Diagram

The 3-color algebra satisfies the following rules:

I. xx = x for all x.
II. (xy)y = x for all x, y.
III. (xy)z = (xz)(yz) for all x,y,z.

As Figure 28 illustrates, these rules correspond to the three Reidemeister moves. An algebra that satisfies these rules is called a quandle (involutory quandle to be precise) [4]. The fact that we can color the trefoil knot with three colors in the 3-color algebra proves that the trefoil knot is knotted. One way to see this is to realize that the quandle rules show that whenever you do a Reidemeister move on a colored diagram, the resulting diagram inherits a unique coloring from the first diagram. From this it is easy to see that any diagram of the trefoil knot that is obtained from the standard diagram has a coloring with three distinct colors. Thus there is no way to transform the trefoil to an unknotted circle by Reidemeister moves, since the circle can be colored with only one color.

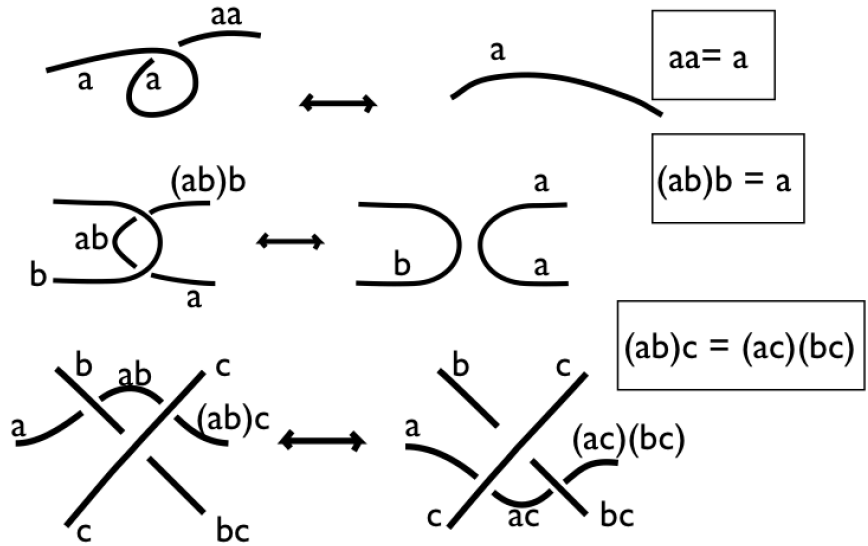

**Figure 28 – The Quandle Rules**

In Figures 29,30 and 31 we illustrate the process of inducing colorings on more complex diagrams of the trefoil knot from an initial coloring of its standard diagram.

# Three-Coloring a Knot

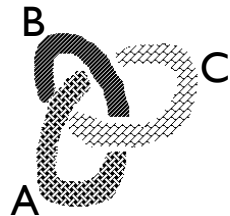

The Rules:
Either three colors at a crossing,
OR
one color at a crossing.

**Figure 29 – Three Coloring Rules for Knots and Links**

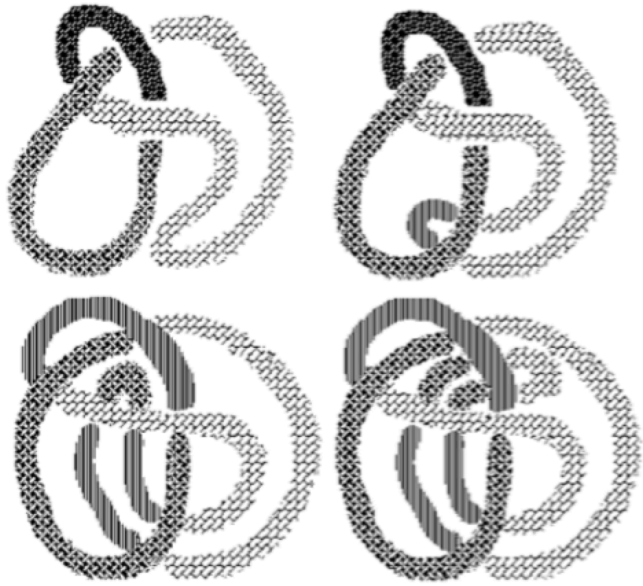

**Figure 30 – Moving the Diagram and Inheriting Color**

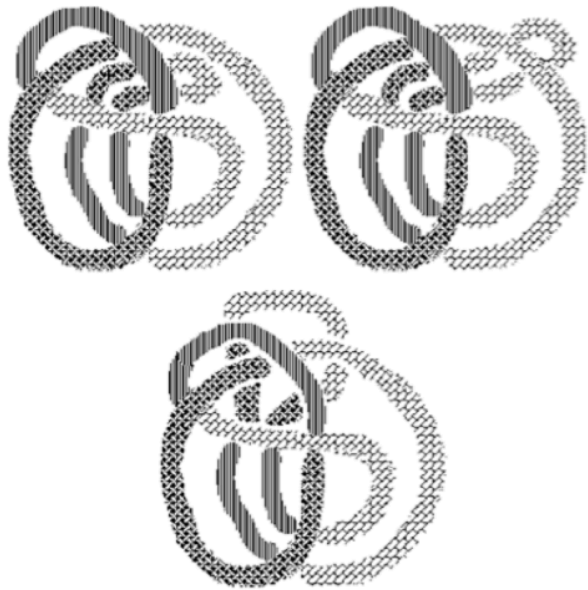

**Figure 31 – Further Moving of the Diagram and Inheriting Colors**

We finish by returning to the Borommean rings of Figures 15, 16 ,17. We can use the 3-color algebra to prove that the Rings are linked. Suppse the Rings are unlinked so that there is a sequence of Reidemeister moves that transforms them to three disjoint circles. Now three disjoint circles can be 3-colored. Apply one color to each circle. Then reverse the sequence of Reidemeister moves and the Borommean rings will inherit a 3-coloring. *If the Rings are unlinked, then they can be non-trivially colored with three colors!*

We leave it as an exercise for the reader to determine that, in fact, *the Borommean rings cannot be non-trivially 3-colored.* And therefore, by the argument given above, they must be linked! This proof of the linkedness of the Borommean rings is due to Colin Adams and will be found in his book[7]. The reader who shows that the Borommean rings can not be colored with three colors can think about how this non-colorability is related to the cyclic surrounding relations of the rings that we discussed in Section 5.

### VIII. A Goedelian Epilogue

We have concentrated mostly on topology and knot theory in this story of mathematical connections.   I would like to end with an example that is purely diagrammatic.  It is well-known that Goedel's incompleteness theorem is based on a method of coding expressions in a formal system by natural numbers in such a way that one can create statements that refer to themselves by referring to their own code numbers. I wanted to diagram the Goedelian situation to clarify it for myself. Accordingly, I chose [8] an arrow of reference in the most general sense.

$$A \rightarrow B$$

The arrow will mean that "A refers to B".  One way to refer is to name, and so the arrow can be interpreted to mean that "A is the name of B." Thus if g is a Goedel number, then

$$g \rightarrow F(u)$$

can mean that "g is the Goedel number of the formula F(u)".  Now with the help of the arrow, we can diagram the famous shift that Goedel devised. Let us assume that the formula F(u) has a single free variable u. (It might be something like "u is a prime number".) Goedel invented a function that I shall denote by #g.  This function is described by the statement "#g is the Goedel number of the formula obtained from F(u) by substituting the Goedel number of F(u) into the free variable of F(u).". That mouthful becomes the following arrow diagram:

$$g \rightarrow F(u), \text{then}$$
$$\#g \rightarrow F(g).$$

Let us call the second arrow the *Goedelian shift* of the first arrow.
Now comes the amazing construction of Goedel. Suppose we have a formula that uses the function #u.  It can have the form F(#u) and it has Goedel number g.

So we have
$$g \rightarrow F(\#u)$$
shifting to
$$\#g \rightarrow F(\#g).$$

Mirablile dictum, the formula F(#g) is discussing its own Goedel number!
This is the key to the Goedelian construction of self-reference. It is the heart of the incompletenesss theorem whereby Goedel creates a sentence that asserts its own unprovability in the given formal system.

The arrow for reference makes the logic of this profound construction easy to survey, and it makes it possible to put the construction in a more general context to see that Goedelian self-reference is very like what we do in language with naming, and eventually naming ourselves. I will stop here.

There are many instances of making connections of ideas and logical relationships by choosing or inventing the right connectors. I hope that the examples in this paper will inspire the reader to collect and create his or her own geometric, topological and diagrammatic connections of mathematical ideas.

**References**
[1] H. Whitney, *On regular closed curves in the plane*. Compositio Mathematica, 4 (1937), p. 276–284.

[2] K. Reidemeister, "Knotentheorie", Chelsea Publishing Co. (1948) (reprint of original 1932 edition).

[3] N. Cozzarelli, J. Dungan, S. Wasserman, Discovery of a predicted DNA knot substantiates a model for site-specific recombination, Science, Vol. 229, (1985), p. 171-174.

[4] L. Kauffman, "Knots and Physics", Fourth edition (2012), World Scientific Pub. Co. Singapore.

[5] J. R. Goldman, L. H. Kauffman, Knots tangles and electrical networks, Advances in Applied Mathematics, Vol. 14, (1993), p. 267-306.

[6] B. Fuller, "Synergetics – Explorations in the Geometry of Thinking", 5.06.00 Knot, p. 229-234. Macmillan Pub. Co. Inc., New York (1982).

[7] C. Adams, "The Knot Book – An Elementary Introduction to the Mathematical Theory of Knots", W. H. Freeman and Co., New York (1994).

[8] L. H. Kauffman. Categorical pairs and the indicative shift. Applied Mathematics and Computation, 218:7989-8004, 2012.